\documentclass[preprint,12pt]{elsarticle}
\usepackage{graphicx}
\usepackage{dsfont}
\usepackage{amsmath,amssymb,amsthm,xcolor,mathabx,float,multirow}
\usepackage{hyperref}
\usepackage{amsbsy}

\usepackage{algorithmic}
\usepackage[ruled,vlined]{algorithm2e}
\usepackage[capitalize,nameinlink,noabbrev]{cleveref}
\usepackage{tikzit}

\tikzstyle{new style 0}=[fill=white, draw=black, shape=circle]
\tikzstyle{new style 1}=[fill=black, draw=black, shape=circle]
\tikzstyle{new style 2}=[fill=none, draw=none, shape=circle]
\tikzstyle{new style 3}=[fill={rgb,255: red,191; green,0; blue,64}, draw=black, shape=circle]

\tikzstyle{new edge style 0}=[-, line width=1 pt]
\tikzstyle{new edge style 1}=[->]
\tikzstyle{new edge style 2}=[-, draw=red, line width=1.2]
\tikzstyle{new edge style 3}=[-, draw=blue, line width=1.2 pt]

\newcommand{\trans  }{^\top}

\newcommand{\Z}{\mathcal{Z}}

\newtheorem{theorem}{Theorem}[section]

\newtheorem{lemma}[theorem]{Lemma}
\newtheorem{Example}[theorem]{Example}
\newtheorem{corollary}[theorem]{Corollary}
\newtheorem{observation}[theorem]{Observation}

\journal{Linear Algebra and Its Applications}

\begin{document}

\begin{frontmatter}

\title{Patterns that Require Distinct Singular Values}


\author[label1]{Caleb Cheung} 
\author[label1]{Bryan Shader}
\affiliation[label1]{organization={Department of Mathematics and Statistics, University of Wyoming},
            addressline={1000 E. University Ave.}, 
            city={Laramie},
            postcode={$82072$}, 
            state={WY},
            country={USA}
            \\
            ccheung@uwyo.edu \quad bshader@uwyo.edu}
\begin{abstract}


\end{abstract}


\end{frontmatter}


\section{Introduction}

Identifying the properties of a structured matrix (e.g., with a given zero-nonzero pattern, or a given are sign-pattern) that allow a multiple eigenvalue is of interest in mechanical engineering (where the presence of a multiple eigenvalue of the stiffness matrix 
allow for more complex vibration patterns), 
and in functional biology (where the presence of a multiple eigenvalue corresponds to multiple equilibria and the ability for a cell to specialize). See 
\cite{HLS} and references therein for more details.

Perhaps the first result along these lines is that of Fiedler \cite{F1964}, which shows that the graphs that forbid a multiple eigenvalue 
are the paths (see below for more specific details).  In this paper, we study the analogous problem for singular values. In particular, we determine the zero-nonzero patterns that forbid a multiple singular value (under a minor natural assumption
on the patterns). 

Throughout the paper,  all matrices are over the reals. The  \textit{singular values} of
the $m\times n$ matrix $A$ are the nonnegative square roots of the eigenvalues of $AA\trans$.
We denote the  multi-set of the singular values of $A$ by $\Sigma(A) = \{\sigma_1,\dots,\sigma_m\}$,
and, as is standard, assume $\sigma_1 \geq \sigma_2 \geq \dots \geq \sigma_m$. 
We  use $O$ to denote a zero matrix and $I$ to denote the identity matrix of the appropriate size. All vectors are  boldfaced, and  $\mathbf{0}$ denotes the zero vector. The \textit{direct sum} of 
the  $m\times n$ matrix $A$ and the  $k \times \ell$ matrix $B$ is denoted by  $A \oplus B$ and is the $(m + k) \times (n+ \ell)$ block diagonal matrix
\begin{align*}
    A \oplus B = \begin{bmatrix}
    \begin{array}{r|r}
        A & O\\ \hline 
        O & B
        \end{array}
    \end{bmatrix}.
\end{align*} 
A \textit{pattern} is a matrix $P=[p_{ij}]$ each of whose entries belongs to $\{0,1\}$.
A \textit{superpattern} of $P$ is 
a pattern obtained by changing some zeros of $P$  to ones. 
The \textit{pattern} of the matrix $A$ is the matrix obtained from $A$ by replacing 
each nonzero entry with a $1$.  The set of all matrices with  pattern $P$ is 
denoted by $\Z(P)$.
 The \textit{term-rank} of the $m\times n$ matrix $A$ is the largest number of 
nonzero entries of $A$ with no two in the same row and column. We say $A$ has \textit{full term-rank} if its term-rank is $\min\{m,n\}$. 

In \cite{FIEDLER1969191}, Fiedler showed that paths are precisely the graphs $G$ for which 
every symmetric matrix whose graph is $G$
has only simple eigenvalues.  In this paper, we study the analogous problem 
for singular values.
In particular, we characterize the $m\times n$ patterns $P$ with $m\leq n$ and term-rank $m$ for which every matrix with pattern $P$ has only simple singular values. 

Unfortunately, allowing a matrix with 
a certain  list of singular values is not inherited by superpatterns. There is a property, which we make use of, that 
 guaranties inheritance. 
The matrix $A$ has the \textit{Strong Singular Value Property} (SSVP) if $O$ is the only matrix $X$ such that
\begin{itemize}
    \item[(a)] $A\trans X$ is symmetric,
    \item[(b)] $X A\trans$ is symmetric, and
    \item[(c)] $A \circ X = O$.
\end{itemize}
Here $\circ$ denotes the entrywise product of matrices. 
Implications of the SSVP are developed 
in \cite{cheung2025strongsingularvalueproperty}.
In particular, 
Theorem 3.10 of that paper shows that every matrix  with the SSVP has full term-rank;
and 
Theorem 4.1  asserts that 
if $A$ has pattern $P$ and the SSVP, then 
for each superpattern $\widehat{P}$ of $P$ there exists a $B$ matrix in $\mathcal{Z}(\widehat{P})$ 
such that $\Sigma(B)=\Sigma(A)$. 

We note that there are related, but fundamentally different investigations in the recent literature, e.g., 
the Strong Inner Product Property \cite{CS2020}, which studies the relationship between patterns and row orthogonality, and the Asymmetric Strong Arnold Property \cite{ARAV2024254} which concerns the combinatorial structure of rectangular matrices and their rank.

The remainder of the paper is structured as follows. In Section 2, we construct families of square patterns that allow
a matrix with multiple singular values  and the SSVP. Using this list, we classify the  square patterns with full term-rank that require all singular values to be simple. 
In Section 3, we extend the characterization to rectangular matrices of full term-rank. In Section 4, we discuss square patterns that allow multiple singular values with less than  full term-rank.

\section{Square Patterns Requiring Simple Singular Values}

 Inverse singular values problems for patterns (ISVP-P) were introduced and studied
in \cite{cheung2025strongsingularvalueproperty}.
Here we recall  needed results from that paper. 
First, it is shown that for each matrix $A$ with the SSVP, and each superpattern 
of that of $A$ there exists  a matrix with the chosen superpattern and the same singular values
as $A$. The formal statement is as follows
(Theorem $4.1$ of
\cite{cheung2025strongsingularvalueproperty}).

\begin{theorem}[\bf Superpattern Theorem]

\label{superpattern}
\phantom{ } \\
    Let $A$ have the SSVP and $P$ be a superpattern of the pattern of $A$. Then there exists a matrix
    with pattern $P$ and singular values $\Sigma(A)$. 
\end{theorem}

Such flexibility is dependent on $A$ having the SSVP. A weaker result, the Matrix Liberation theorem, allows one to perturb a matrix $A$ without the SSVP in certain directions to a matrix $B$ with  
$\Sigma(B) = \Sigma(A)$.
To understand the statement of the theorem, we introduce some definitions.

We view the set $\mathbb{R}^{m\times n}$ of all $m\times n$ real matrices as an inner-product space with inner-product $\langle A, B \rangle= \mbox{tr}(A\trans B)$, where $\trans$ and $\mbox{tr}$ indicate the transpose and trace of a matrix, respectively. Let $A$ be an $m\times n$ matrix, and $S$ be an $m\times n$ superpattern of the pattern of $A$. The matrix $A$ has the {\it SSVP with respect to $S$} provided $Y=O$ is the only $m\times n$ matrix such that 
\begin{itemize}
\item[(a)]
$A\trans Y$ is symmetric,
\item[(b)]
$YA\trans$
is symmetric, and 
\item[(c)] $S \circ Y=O$.
\end{itemize} 
The \textit{tangent space of $A$}, denoted $\mbox{Tan}^{\Sigma}_A$, is the collection of matrices of the form 
$KA+AL$ where $K$ is an $m\times m$, respectively $L$ is an $n\times n$, skew-symmetric matrix. The \textit{normal space} of $A$, denoted $\mbox{Norm}^{\Sigma}_A$ is the orthogonal complement of the tangent space of $A$. In 
\cite{cheung2025strongsingularvalueproperty} it is shown that 
$\mbox{Norm}^{\Sigma}_A$ is the set of all matrices $X$ such that 
both $A\trans X$ and $XA\trans$ are symmetric. For an $m\times n$ pattern $P=[p_{ij}]$,  and an $m\times n$ matrix
$D=[d_{ij}]$,
the {\it pattern of $P$ in the direction $D$} is the $m\times n$ 
pattern $S=[s_{ij}]$ with $s_{ij}=1$ if and only if  $p_{ij}=1$ or $d_{ij}\neq 0$.
The following is Theorem 4.5 of 
\cite{cheung2025strongsingularvalueproperty}

\begin{theorem}[\bf Matrix Liberation Theorem]
\phantom{ } \\
Let $A$ be an $m\times n$ matrix with pattern $P$,  $D$ be a matrix in $\mbox{\rm Tan}^{\Sigma}_A$ 
and  $S$ be the pattern of $P$ in the direction of $D$.
If $A$ has the SSVP with respect to $S$, then there is a matrix $B$ with pattern 
$S$ having the SSVP, the  same singular values as $A$ and pattern $S$.

\end{theorem}

Another tool we need is the Direct Sum theorem, which classifies the matrix direct sums having the SSVP (Theorem 3.12 of 
\cite{cheung2025strongsingularvalueproperty}).

\begin{theorem}[\bf Direct Sum Theorem]
\label{dsum}
\phantom{ }\\
Let $A$ be an $m\times n$ matrix, $B$ be a $p \times q$ matrix and $M=A\oplus B$ with $m+p \leq n+q$.
Then $M$ has the SSVP if and only if  each of the following holds
\begin{itemize}
\item[\rm (a)] 
$ m \leq n$ and $p\leq q$;
\item[\rm (b)] 
both $A$ and $B$ have the SSVP; 
\item[\rm (c)] $A$ and $B$ have no common nonzero singular value; and
\item[\rm (d)] either both $A$ and $B$
have linearly independent rows, or one of $A$ and $B$ is square and invertible.
\end{itemize}
\end{theorem}

A necessary condition for an $m\times n$ matrix $A$ with $m\leq n$ to have the SSVP is that $A$ have term-rank $m$ (see Theorem 3.10 of \cite{cheung2025strongsingularvalueproperty}).  Thus, we primarily  restrict our study to matrices with full term-rank.

Let $M$ be an $n\times n$ matrix with full term-rank. We say that $M$ is in {\it standard form} provided each of its diagonal entries is nonzero. Note there is no loss of generality 
in assuming $M$ is in standard form, since 
any $n$ nonzero entries of $M$ in distinct rows and columns can be placed on the diagonal by permuting rows of $M$, 
and such an operation preserves singular values and the SSVP. When $M=[m_{ij}]$ is in standard form, it is convenient to associate $M$ with 
its digraph; the  digraph having vertices $1$, \ldots, $n$ and an arc from $i$ to $j$ 
if and only if $i\neq j$ and $m_{ij}\neq 0$.
Note in describing the digraph of $M$ 
diagonal entries are ignored; since each is nonzero.

We now focus on identifying those $n\times n$
patterns with full term-rank that do not allow a matrix with a multiple singular value.  The plan of attack is to find matrices in standard form with the SSVP and a multiple singular value by direct computation or via the Matrix Liberation theorem. Then, by the Direct Sum and Superpattern theorems, the patterns of such matrices are  forbidden  principal submatrices for patterns in standard form  that require all singular values to be simple. Then we  characterize all patterns $Q$ (in standard form) that avoid these  forbidden patterns and argue that each matrix with such a pattern has only simple singular values.

We need the following lemma that describes the solutions to Sylvester's equation $AX=XB$
in the case that $A$ and $B$ are symmetric $m\times m$ and $n\times n$ matrices, respectively. 
\begin{lemma}
\label{sylvester}
Let $A$ and $B$ be symmetric $m\times m$ and $n\times n$ 
matrices, respectively, with spectral decompositions
$A=\sum_{i=1}^m \lambda_i \mathbf u_i \mathbf u_i\trans$, 
and 
$B=\sum_{j=1}^n \mu_j\mathbf v_i \mathbf v_i\trans$.  Then the $m\times n$ matrix $X$ 
satisfies $AX=XB$ if and only if there exist scalars
$c_{ij}$ such that $X=\sum c_{ij} \mathbf u_i \mathbf v_j\trans$, where the sum is over all $(i,j)$ where 
$\lambda_i=\mu_j$.
\end{lemma}

\begin{proof} 
Observe that $\mathbf u_i\mathbf v_j\trans$ $(i=1,\ldots, m, j=1,\ldots, n)$
is a basis of $\mathbb{R}^{m\times n}$.  Let $X= \sum_{i=1}^m\sum_{j=1}^n
c_{ij} \mathbf u_i \mathbf v_j\trans$ be an arbitrary matrix in $\mathbb{R}^{m\times n}$.
Then 
\begin{align}
AX&= \sum_{i=1}^n \sum_{j=1}^n
\lambda_i c_{ij}  \mathbf u_i \mathbf v_j\trans \mbox{ and}\\
XB&= \sum_{i=1}^n \sum_{j=1}^n \mu_j c_{ij}  \mathbf u_i \mathbf v_j\trans. 
\end{align} 
Thus, $AX=XB$ if and only if $c_{ij}=0$ 
whenever $\lambda_i\neq \mu_j$. 
\end{proof}

A square matrix $N$ is \textit{irreducible} if it is not permutationally similar to a matrix of the form 
\[ 
\begin{bmatrix}
\begin{array}{c|c}
    N_{1} & N_{12} \\ \hline
    O & N_{2}
    \end{array}
\end{bmatrix},
\]
where $N_1$ and $N_2$ are square matrices each of order at least $1$.  The Perron-Frobenius theorem 
implies that if $N$ is an entrywise nonnegative 
irreducible matrix, then it has an entrywise positive 
eigenvector corresponding to its spectral radius.

The {\it bigraph} of the $m\times n$ pattern $P$ is  denoted by 
$\mathcal{B}(P)$, has vertex set $\{1,2,\ldots, m,1',2', \ldots, n'\}$, and an edge joining vertex $i$ to vertex $j'$ if and only if $a_{ij}\neq 0$. We call $1, \ldots, m$
the \textit{row vertices} of $\mathcal{B}(P)$ 
and $1', \ldots, n'$ the \textit{column vertices} of $\mathcal{B}(P)$. The \textit{row graph} of $P$
consists of vertices $1, \ldots, m$ with an 
edge joining $i\neq j$ if and only if there is $k$ such that $p_{ik}=1=p_{jk}$.  Evidently, if $P$ has 
no row of zeros and no column of zeros, then the bigraph of $P$ is connected if and only if the 
row graph of $P$ is connected; and hence if the bigraph of $P$ is connected and $A$ is 
an entrywise nonnegative matrix  $\mathcal{Z}(P)$, 
then $AA\trans$ is an irreducible matrix.

\begin{theorem}
\label{sumplus}
Let $P$ and $Q$ be full term-rank matrices
of size $m\times m$ and $n\times n$ respectively 
such that 
each of their bipartite graphs is connected.
  Let $E$, respectively $F$, be an $m\times n$, 
respectively $n\times m$  pattern 
such that the total number of nonzero entries in $E$ or $F$ is at least $2$. 
Then there is a matrix with pattern 
\begin{equation}
\label{patt}
M=\begin{bmatrix}
\begin{array}{c|c}
P & E\\ \hline
F & Q
\end{array}
\end{bmatrix}
\end{equation}
having the SSVP and a multiple singular value. 
\end{theorem}

\begin{proof}
Let $1=\sigma_1> \cdots > \sigma_m>0$
and $1=\tau_1> \cdots > \tau_n>0$ 
be real numbers such that $\sigma_i=\tau_j$ if and only if $i=j=1$. 
By the Direct Sum theorem, and the Supergraph theorem 
there exists a matrix $A \in \mathcal{Z}(P)$ 
with singular values $\{ \sigma_1, \ldots, \sigma_m\}$.
Indeed, the ``signed'' version of the Supergraph 
theorem, guarantees that there is such an $A$ 
that is entrywise nonnegative matrix.  Similarly, there exists an entrywise nonnegative matrix $B 
\in \mathcal{Z}(Q)$ with singular values 
$\tau_1, \ldots, \tau_n$.  Since each column and row of  $P$ and $Q$ is nonzero and the bigraphs of $P$ 
and $Q$ are connected, the eigenspace of each of  $AA\trans$, $A\trans A$, $BB\trans$ and $B\trans B$ 
corresponding to $1$ is spanned by an entrywise positive vector. 

We wish to apply the Matrix Liberation theorem to $A \oplus B$.  We begin by determining $\mbox{Tan}_{A\oplus B}^{\Sigma}$.  We do this by determining its orthogonal complement  $\mbox{Norm}_{A \oplus B}^{\Sigma}$ in the set  of $(m+n) \times (m+n)$ matrices.
We know that 
$\mbox{Norm}^{\Sigma}_{A \oplus B}$ is the 
set of 
all matrices of the form 
\[
X= 
\begin{bmatrix}
\begin{array}{c|c}
R& S \\  \hline 
T & U 
\end{array}
\end{bmatrix},
\]
where $R$ is $m \times m$, $U$ is $n\times n$,
 $(A \oplus B)\circ X=O$, 
and both $(A\oplus B)\trans X$ and $X \trans (A\oplus B)$ are symmetric. Consider such an $X$.
Then $A\circ R=O$, and both $A\trans R$
and $R A\trans$ are symmetric. We conclude that $R=O$.  Similarly, $U=O$.  Additionally, 
we have that 
\begin{align}
\label{one}
S\trans A&=B\trans T, \mbox{ and }\\ 
\label{two}
BS\trans &= T A \trans .
\end{align}
Pre-multiplying (\ref{one}) by $B$ and substituting in 
(\ref{two}) gives
\[ BB\trans T=BS\trans A=TA\trans A.\]
By \Cref{sylvester}, $T= c\mathbf u \mathbf v\trans$ where $\mathbf u$ is an eigenvector of $BB\trans$ corresponding to the eigenvalue $1$,  $\mathbf v$ is 
an eigenvector of $A\trans A$ corresponding to the eigenvector $1$, and $c$ is a constant. 
Similarly, $S\trans=d \mathbf w \mathbf x\trans$, 
where $\mathbf x$ is an eigenvector of $A A\trans$ corresponding to $1$, $\mathbf w$
is an eigenvector of $B\trans B$ corresponding to $1$, and $d$ is a constant. 
As  noted above, $\mathbf u$, $\mathbf v$, $\mathbf w$ and $\mathbf x$
can be taken to be entrywise positive. 

Since $B$ 
is invertible (\ref{one}) implies that 
$S=(A^{-1})\trans T\trans B$.  In particular, 
we have shown that $\mbox{Norm}^{\Sigma}_{A \oplus B}$ is contained in the span of a single matrix of the form 
\begin{equation}
\label{spanner}
\begin{bmatrix}
\begin{array}{cc}
O & W
\\
Z& O 
\end{array}
\end{bmatrix},
\end{equation}
where $W$ and $Z$ have no zero entries.
Since there are at least two nonzero entries 
in $E$ and $F$, the one-dimensionality of 
$\mbox{Norm}^{\Sigma}_{A \oplus B}$ and 
the fact that each entry of $W$ and $Z$ is nonzero imply  that there is a matrix with pattern 
\begin{equation}
\label{direct}
\begin{bmatrix}
\begin{array}{cc}
O & E\\
F & O
\end{array}
\end{bmatrix}
\end{equation}
that is orthogonal to the matrix in (\ref{spanner}).  Such a matrix is in
$\mbox{Tan}_{A\oplus B}^{\Sigma}$. 
Additionally, if one includes the requirement 
that 
$S \circ E=O$ and $T \circ F$, then the above analysis 
implies that  $A \oplus B$ has the  strong property with respect to the graph of $A \oplus B$ 
in the direction of the matrix in (\ref{direct}). 
The Matrix Liberation theorem now implies that 
there is a matrix with pattern $M$ having the SSVP and the same singular values as $A \oplus B$.  
\end{proof}

\Cref{sumplus} has some useful consequences 
for $n\times n$ patterns $P$  with term-rank $n$ 
in standard form. A \textit{weak cycle} of the digraph of $P$ is a collection of distinct arcs in the digraph
such that ignoring the direction of the arcs of the results in a graph that is a cycle (possibly of length 2).

\begin{corollary}
\label{cycleforbid}
Let $W$ be an $n\times n$ pattern in standard
form whose digraph contains a weak cycle. 
Then there exists a matrix in $\mathcal{Z}(W)$
with the SSVP and a multiple singular value.
\end{corollary}

\begin{proof}
By the Direct Sum theorem and the Superpattern theorem,
it suffices to prove the result when $W$ is in standard form and its digraph is a weak cycle. 
Note that if we replace the two off-diagonal 1's
in  the first row or column of $W[\alpha]$ by zeros, 
then the resulting pattern, $[1] \oplus W[\alpha\setminus \{1\}]$,
satisfies the hypothesis of \Cref{sumplus}, and the result follows. 
\end{proof}

In addition to those patterns in \Cref{cycleforbid},
there are other normalized square patterns 
that allow a matrix with the SSVP and a multiple singular.  We describe two such patterns in the next lemmas.

\begin{lemma}
\label{deg4}
   The  pattern
    \[ P=
    \left[ 
    \begin{array}{cccc}
         1 & 1 & 1 & 1  \\
         0& 1 & 0 & 0\\
         0 & 0 & 1 & 0\\
         0 & 0 & 0 & 1\\
    \end{array}
    \right]   
    \]
allows  a matrix with multiple singular values and the SSVP.
\end{lemma}
\begin{proof}
Note that
     \[ PP\trans- I=
    \left[ 
    \begin{array}{cccc}
         3 & 1 & 1 & 1  \\
         1& 0 & 0 & 0\\
         1 & 0 & 0 & 0\\
         1 & 0 & 0 & 0
    \end{array}
    \right]   
    \]
has nullity $2$. 
Hence $P$ has $1$ as a singular value of multiplicity 2. Now let
     \[ X=
    \left[ 
    \begin{array}{cccc}
         0 & 0 & 0 & 0  \\
         a& 0& b & c\\
         d & e& 0& f\\
         g & h & i & 0\\
    \end{array}
    \right]   
    \]
be such that $XP\trans$ and $P\trans X$ are both symmetric. 
The latter gives
      \begin{align*}
      a,d,g &= 0 \\
      b & =  e \\
      c & =  h\\
      f & =  i.\\
\end{align*}
Additionally, the former gives
\begin{align*}
    a+b+c & =0\\
    d+e+f & = 0\\
    g+h +i &= 0.
\end{align*}
Together these imply,
\begin{align*}
\begin{bmatrix}
\begin{array}{ccc}
1 & 1 &0 \\
1 & 0 & 1\\
0 & 1 & 1 \end{array}
\end{bmatrix}
\begin{bmatrix}
\begin{array}{c}
b\\
c\\
f
\end{array} 
\end{bmatrix}= \begin{bmatrix}
    0\\
    0\\
    0
\end{bmatrix}.
\end{align*}
Hence $b=c=f=0$, which implies $e=h=i=0$.
Thus $X = O$, and $P$ has the SSVP.
\end{proof}

\begin{lemma}
\rm
\label{inout}
The pattern
\[ A=\left[\begin{array}{rrrr}
1 & 1 & 0 & 0 \\
0 & 1 & 1 & 1 \\
0 & 0 & 1 & 0 \\
0 & 0 & 0 & 1
\end{array}\right]
\]
allows a matrix that has the SSVP and a multiple singular value. 
\end{lemma}

\begin{proof}
Consider the matrix
\[ A=\left[\begin{array}{rrrr}
1 & 1 & 0 & 0 \\
0 & 1 & 1 & 1 \\
0 & 0 & \sqrt{2} & 0 \\
0 & 0 & 0 & \sqrt{2}
\end{array}\right].
\] 
Since $(AA\trans)(\{2\})=2I_3$, $\sqrt{2}$ is a multiple singular value of $A$. To show that $A$ has the SSVP, consider a matrix of the form 
\[ 
X=\left[\begin{array}{rrrr}
0 & 0 & a & b \\
c & 0 & 0 & 0 \\
d & e & 0 & f \\
g & h & i & 0
\end{array}\right]
\]
for which both $AX\trans$ and $X\trans A$ are symmetric. 
These yield the following system of equations
\[ 
\begin{array}{lllll}
c = 0,&   c = a+b, \\
 a= \sqrt{2}d+c, & \sqrt{2}a=d+e,\\
 b=\sqrt{2}g+c, & \sqrt{2}b=g+h\\
 a= \sqrt{2}e, & e+f=0\\
 b= \sqrt{2}h, & h+i=0, \\ 
 \sqrt{2}i=\sqrt{2}f. \\
\end{array}
\] 
whose only solution is the trivial one.  Hence $A$ has the SSVP.
\end{proof}

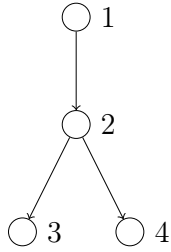
\begin{figure}
\begin{center} 
  \resizebox{1 in}{!}{
\begin{tikzpicture}[scale=.75]
	\begin{pgfonlayer}{nodelayer}
		\node [style=new style 0, minimum size=.05cm, label=right:2] (2) at (0, 1) {};
		\node [style=new style 0, minimum size=.05cm, label=right:3] (3) at (-1, -1) {};
		\node [style=new style 0, minimum size=.05cm,label=right:4] (4) at (1, -1) {};
		\node [style=new style 0, minimum size=.05cm, label=right:1] (5) at (0, 3) {};
	\end{pgfonlayer}
	\begin{pgfonlayer}{edgelayer}
		\draw [style=new edge style 1] (2) to (3);
		\draw [style=new edge style 1] (2) to (4);
		\draw [style=new edge style 1] (5) to (2);
	\end{pgfonlayer}
\end{tikzpicture}
}
\caption{
  {Digraph of pattern in \Cref{inout}}}
    \label{fig1}
\end{center}
\end{figure}

The digraph of the pattern in \Cref{inout} is
given in \Cref{fig1}.

\begin{corollary}
\label{degforbid}
If $P$ is an $n\times n$ pattern in standard form, 
and the digraph of $D$ contains a vertex whose 
in-degree and out-degree sum to at least $3$, 
then there is a matrix in $\mathcal{Z}(P)$
with a multiple singular value. 
\end{corollary}
\begin{proof}
 This follows from the  Direct Sum theorem,  Superpattern theorem, \Cref{deg4}, and  \Cref{inout}
 \end{proof}

A digraph is a \textit{weak path} provided 
its underlying graph is a path. We now show that for $n\times n $ patterns $P$ with term-rank $n$ a necessary condition to require all simple singular values is that the digraph of $P$ is a weak path.

\begin{theorem}
\label{th:weakpath}
Let $P$ be an $n\times n$ pattern with term-rank $n$ in standard form.
If $P$ does not allow a matrix with a multiple singular value, then the digraph of $P$ 
is a weak path. 
\end{theorem}

\begin{proof}
Assume $P$ does not allow a matrix with a multiple singular value.
Let $\mathcal{D}$ be the digraph of $P$, and $G$ be the underlying graph of $\mathcal{D}$.  

Necessarily, $G$ is connected; for otherwise, 
$P$ is permutationally similar to a direct sum of patterns in standard form and one can make each of the direct summands share their largest singular value to obtain a matrix with pattern $P$ and a multiple singular value.

By \Cref{degforbid} and \Cref{cycleforbid},
each vertex of $\mathcal{D}$  has the sum of its in-
and out-degree at most two, and $\mathcal{D}$
has no weak-cycle. This implies that 
$\mathcal{D}$
is a weak path.
\end{proof}

Our next task is to show that each pattern $P$ in standard form 
whose digraph is a weak path
does not allow a multiple singular value.
If the digraph of $P$ is a directed path, then for every matrix $A$
with pattern $P$, the graph of the matrix $AA\trans $ is a path, 
and hence its eigenvalues (and the singular values of $A$) are 
distinct.  The next example shows that a pattern $P$ can require all simple singular values even when the graph of $PP\trans$ allows a multiple eigenvalue.

\begin{Example}
\rm 
Consider the pattern 
\[ 
P=\left[ \begin{array}{cccc}
1 & 1 & 0 & 0\\
0 & 1 & 0 & 0\\
0 & 1 & 1 & 1\\
0 & 0 & 0 & 1
\end{array}
\right].
\]

We claim that each matrix with this pattern has distinct singular values.

 Let $a,b,c,d,e,f,g$ are nonzero real numbers
such that 
\[ 
M=\left[ \begin{array}{cccc}
a & b & 0 & 0\\
0 & c & 0 & 0\\
0 & d & e & f\\
0 & 0 & 0 & g
\end{array}
\right]
\]
and let $\sigma$ be a singular value of $M$.  
Note that $M$ is nonsingular, so $\sigma \neq 0$. 
Also, 
\[ 
MM\trans-\sigma^2 I=
\left[ \begin{array}{cccc}
a^2+b^2-\sigma^2 & bc & bd & 0\\
bc & c^2-\sigma^2 & cd & 0\\
bd & cd & d^2+e^2 + f^2 -\sigma^2 & fg\\
0 & 0 & fg & g^2-\sigma^2
\end{array}
\right].
\]
Now note that the matrix 
\[ 
N=
\left[ \begin{array}{cc|c}
 bc & bd & 0\\
 c^2-\sigma^2 & cd & 0\\ \hline
cd & d^2+e^2 + f^2 -\sigma^2 & fg\\
\end{array}
\right].
\]
obtained from $MM\trans -\sigma^2 I$
by deleting its last row and first column is block lower-triangular. Since $\sigma \neq 0$, the leading $2\times 2$ principal submatrix of $N$  is invertible. 
 Hence $N$ is invertible, $\sigma^2$ is an eigenvalue of $MM\trans$ of multiplicity $1$, and  $\sigma$ is a singular value of $M$ 
of multiplicity $1$. The graph of $MM\trans$ is the Paw. Thus there
is a symmetric (positive semi-definite) matrix $A$ whose graph is the Paw
having a multiple eigenvalue.  Hence $A=BB\trans$ for some $4\times 4$ matrix $B$. The pattern of $B$ is necessarily not that of $P$. 
$\diamond$
\end{Example}

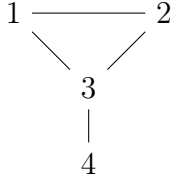
\begin{figure}
\begin{center}
\begin{tikzpicture}
     \node (v1) at (0,0) {1};
     \node (v2) at (1,-2) {4};
     \node (v3) at (2,0) {2};
     \node (v4) at (1,-1) {3};
     \draw (v1) -- (v3) -- (v4)--(v1);
     \draw (v2) -- (v4); 
\end{tikzpicture}
\\
\caption{The Paw Graph}
\label{paw}
\end{center}
\end{figure}

Let $D$ be a weak path on  $n$ vertices. 
Then there exist $1=i_0 <i_1< \cdots <i_k=n$
such that between $i_{j-1}$ and $i_j$ there 
is either a directed path $P_j$ from $i_{j-1}$ 
to $i_j$ or a directed path $P_j$ from $i_{j}$ to $i_{j+1}$ in $D$ for $j=0, \ldots, k-1$ and each $i_j$ is either has in-degree zero or 
out-degree zero $j=0,\ldots, k$.

\begin{center}
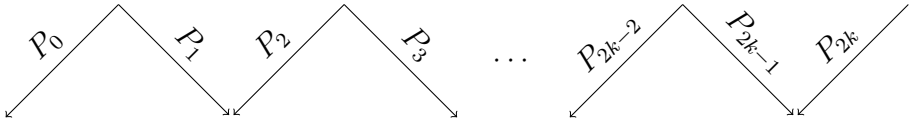
\begin{figure}
\begin{tikzpicture}[scale=1.5]
\draw[->] (0,0)--(-1,-1) node [pos=0.5, above, sloped] {$P_0$};
\draw[->] (0,0)--(1-.02,-1+.02) node [pos=0.5, above, sloped] {$P_1$};
\draw[->] (2,0)--(1+.02,-1+.02) node [pos=0.5, above, sloped] {$P_2$};
\draw[->] (2,0)--(3,-1) node [pos=0.5, above, sloped] {$P_3$};
\draw[->] (5,0)--(4,-1) node [pos=0.5, above, sloped] {$P_{2k-2}$};
\draw[->] (5,0)--(6-.02,-1+.02) node [pos=0.5, above, sloped] {$P_{2k-1}$};
\draw[->] (7,0)--(6+.02,-1+.02) node [pos=0.5, above, sloped] {$P_{2k}$};
\node at (3.5,-.5) {$\cdots$};
\end{tikzpicture}
\\
\caption{
 Decomposition into directed paths}
\label{zigzag}
\end{figure}
\end{center}

Now assume that $P$ is an $n\times n$ pattern 
with term-rank $n$ in standard form whose digraph $D$ is a weak path, 
and that  $D$ is decomposed into
directed paths as described above. 
Note for each $A \in \mathcal{Z}(A)$ the graph of $AA\trans$ 
is that of $PP\trans$. The graph 
of $PP\trans$ has vertices $1,\ldots, n$ 
and an edge joining $i\neq j$ if and only if 
row $i$ of $P$ and row $j$ of $P$ have nonzero 
dot-product. That is, if and only if there is an arc in $D$ between vertex $i$ and vertex $j$,  or there is a vertex $k$ such that both $ij$ and $ik$ are arcs of $D$. Hence, the 
graph of $PP\trans$ consists of the path 
$1$--$2$---$\cdots$---$n$ along with all edges of the form $uv$ where $ui$ and $vi$ are arcs
and $i$ is the terminal vertex of two directed paths in $D$. 

For example, if $D$ is the digraph in 
\Cref{zigzag}, 
then the corresponding graph is in \Cref{row}. 

\begin{center}
\begin{figure}
\begin{tikzpicture}[scale=1.5]
\draw[-] (0,0)--(-1,-1){};
\draw[-] (0,0)--(1,-1) node [pos=0.7, below, sloped] {$v_1$};
\draw[-] (2,0)--(1,-1) node [pos=0.7, below, sloped] {$v_2$};
\draw[-] (2,0)--(3,-1)  ;
\draw[-] (5,0)--(4,-1)  ;
\draw[] (.8,-.8)--(1.2,-.8);
\draw[-] (5,0)--(6,-1) node [pos=0.7, below, sloped] {$v_{2k-1}$};
\draw[-] (7,0)--(6,-1) node [pos=0.7, below, sloped] {$v_{2k}$};
\draw[] (5.8,-.8)--(6.2,-.8);
\node at (3.5,-.5) {$\cdots$};
\node at (1,-1.2){};
\node at (6,-1.2) {};
\end{tikzpicture}
\\
\caption{Graph of $PP\trans$}
\label{row}
\end{figure}
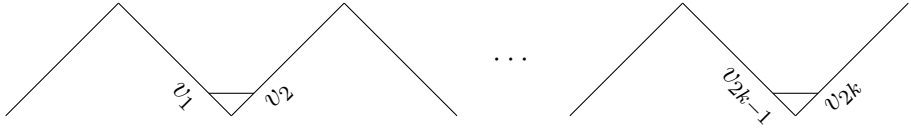

\end{center}

\begin{lemma}
\label{zigzagdistinct}
    Let $P$ be an $n\times n$ pattern with term-rank $n$ in standard form whose digraph 
    is a weak path. Let $A \in \mathcal{Z}(P)$.
    Then the singular values of $A$ are distinct. 
\end{lemma}
\begin{proof}
Without loss of generality, we may assume that the vertices are ordered so that the underlying graph of $P$ 
is the path $1$--$2$---$\cdots$---$n$.

Let $s_1,s_2,\dots,s_k$ be the vertices of the digraph of $D$ of $P$ of outdegree $0$. 
Since the digraph of $P$  is a weak path, the graph of $AA\trans$ consists of the path $1$---$2$---$\cdots$---$n$ along with  an  edge between 
the neighbors of $s_i$ for  $i=1,\ldots, k$.
Evidently, $(AA\trans)(n,1)$ is a block lower-triangular matrix whose diagonal blocks are $1\times 1$
blocks corresponding to superdiagonal entries of $AA\trans$ and $2\times 2$
blocks of the form
\begin{align}
\label{2by2}
\left[ 
\begin{array}{cc}
     a_{s_{i}-1,s_i}a_{s_{i},s_i} & a_{s_{i-1},s_i}a_{s_i+1,s_{i}}  \\
     a_{s_{i},s_i}^2 & a_{s_{i}
     ,s_i}a_{s_i+1,s_{i}}
\end{array}
\right],
\end{align}
for $j=1,\ldots, k$.

Suppose $\lambda$ is an eigenvalue of $AA\trans$. Then the $1\times 1$ diagonal blocks
of $(AA\trans)(n,1)$ are nonzero and hence invertible. Each of the $2 \times 2$ blocks (\ref{2by2}) has the form
\[
\left[ 
\begin{array}{cc}
     bd & de  \\
     b^2-\lambda & be
\end{array}
\right],
\]
which has determinant $de\lambda$. Since $A$ is invertible, $\lambda \neq 0$. Hence each diagonal block of $(AA\trans-\lambda I)(n,1)$ is invertible, 
and the nullity of $AA\trans -\lambda I$ is 1. Therefore, $A$ has no multiple singular value.
\end{proof}

 \Cref{th:weakpath}
and \Cref{zigzagdistinct} give the
following characterization for square 
patterns with full term-rank that require  singular values to be  distinct.

\begin{corollary}
\label{fiedlersquare}
Let $P$ be an $n\times n$ pattern with term-rank $n$ and assume $P$ is in standard form with digraph $D$.  Then every matrix in $\mathcal{Z}(P)$ has all simple singular values
if and only if $D$ is a weak path. 
\end{corollary}

\section{Non-square patterns requiring all simple singular values}

We now turn our attention to $m\times n$ patterns $P$ with $m<n$ having full term-rank
and the property that each singular value of each matrix with pattern $P$ has multiplicity $1$. We begin with  examples that show the subtlety of this problem.

\begin{Example}
\label{forbiddensun}
\rm
Consider the pattern 
\[
P= \begin{bmatrix}
\begin{array}{cccc}
1 & 0 &0  & 1 \\
1 & 1 & 0 & 0\\
1 & 0 & 1 & 0
\end{array}
\end{bmatrix}.
\]
The leading $3\times 3$ pattern $Q$ of $P$
has digraph that is a weak path. It can be verified that each matrix in $\mathcal{Z}(Q)$ 
has the SSVP. However,  by \Cref{zigzagdistinct},  no matrix in $\mathcal{Z}(Q)$ has a multiple singular value. Hence,
it is not possible to use $Q$ and the Superpattern theorem to argue that $P$ allows a matrix with a multiple singular value. 

One can  permute the rows and columns of $P$ to get the matrix $R=[ I_3 \; J_{3,1} ]$.
The leading $3\times 3$ matrix $I_3$ does allow matrices with multiple singular values, but none of these has the SSVP. 
Yet, again
it is not possible to use the Superpattern  theorem to argue that $P$ allows a matrix with a multiple singular value. 

The two 
instances above are essentially the only ways to permute the rows and columns of $P$ to get the leading $3\times 3$ submatrix in standard form.
So, either $P$ does not allow a matrix with 
a multiple singular value, or we need to 
find and verify such a matrix without using
the SSVP.  
Note that  $RR\trans= I+J$, and hence $R$ has 
singular values $2$, $1$ and $1$! While in this setting, we also verify that $R$ has the SSVP. 
To see this, let 
\[ 
X= \begin{bmatrix}
\begin{array}{cccc}
0 & a & b & 0 \\
c & 0 & d & 0 \\
e& f & 0 & 0 
\end{array} 
\end{bmatrix},
\]
and assume that both $R\trans X$ and $XR\trans$
are symmetric. Then
$X[\{1,2,3\}]$ is symmetric with row sums zero.
Evidently, this implies that $X=O$, and 
$R$ has the SSVP. 

\hfill $\diamond$  \end{Example}

\begin{Example}
\label{aster}
\rm 
Let 
\[ P = 
\begin{bmatrix}
\begin{array}{cccc}
1 & 0 & 0 & 0\\
1 & 1 & 0 & 0 \\
1 & 0 & 1 & 1 
\end{array} 
\end{bmatrix}.
\]
Then $P[\{1,2,3\}]$ is in standard form and its digraph is a weak path. Let $A \in \mathcal{Z}(P)$. 
There exists a $2\times 2$ orthogonal matrix 
$Q$ such that the pattern of 
$A(I_2 \oplus Q)$ is that of $P$ with the $(3,4)$-entry replaced by a  $0$.  Hence, the 
singular values of $A$ are equal to the singular values of a $A(I_2 \oplus Q)[\{1,2,3\}]$ whose pattern is a weak path.  Hence the singular values of $A$ are distinct. 

Therefore, every matrix with pattern $P$ 
has each of its singular values of multiplicity $1$. 
 $\diamond$  \end{Example}

More generally, we have the following. 

\begin{lemma}
\label{columnaugment}
    Let $P$ be an $m\times n$ pattern with $m \leq n$ having a column with exactly one $1$. 
    Let $P^+$ be the $m\times (n+1)$ pattern obtained from $P$ by appending a duplicate of that column.
    Then 
    \begin{itemize}
    \item[\rm (a)]  For each matrix  $B \in \mathcal{Z}(P^+)$
     there is a matrix $A \in \mathcal{Z}(P)$ with   $\Sigma(A) = \Sigma(B)$. 
    \item[\rm (b)] For each matrix  $A \in \mathcal{Z}(P)$
     there is a matrix $B \in \mathcal{Z}(P^{+})$ with   $\Sigma(A) = \Sigma(B)$. 
    \item[\rm (c) ] $P$ allows a matrix with a multiple singular value if and only if $P^+$ does. 
     \end{itemize}
\end{lemma}
\begin{proof}
Without loss of generality, the last two columns of $P$ and $P^+$ 
are equal columns of support $1$.

Let $B \in \mathcal{Z}(P^{+})$. The last two columns of $B$ are 
$c \mathbf e_i$ and $d \mathbf e_i$ for some $i$ and some nonzeros $c$ and $d$.
Let 
\[ Q= I_{n-2} \oplus 
\frac{1}{\sqrt{c^2+d^2}}
\left[ \begin{array}{rr}
c& d\\ 
d &   -c
\end{array}
\right]. \]
Then $Q$ is an orthogonal matrix, the 
first $n-1$ columns of $BQ$ are those of $B$,
the $n$-th column  is $\sqrt{c^2+d^2}\mathbf e_i$, and the last column is $\mathbf 0$.
Hence $A=B[\{1,\ldots,m\},\{1,\ldots, n\}]$ has the same singular values as $B$, and (a) holds.

Now let $A \in \mathcal{Z}(P^{+})$. Append a column of zeros to the right of $A$ to get a matrix $A'$.
The matrix $B=A'Q$ where $Q$ is the orthogonal matrix
\[ 
I_{n-2} \oplus 
\frac{1}{\sqrt 2}\begin{bmatrix}
\begin{array}{rr}
    1 & -1 \\
    1 & 1
    \end{array}
\end{bmatrix}
\]
has pattern $P^+$ and the same singular values as $A$.  Hence (b) holds.

Statement (c) follows from (a) and (b). 
\end{proof}

\begin{figure}
\begin{center}
\begin{tikzpicture}[scale=.5]
	\begin{pgfonlayer}{nodelayer}
		\node [style=new style 0] (0) at (0, 3) {};
		\node [style=new style 1] (1) at (0, 1.5) {};
		\node [style=new style 0] (2) at (0, 0) {};
		\node [style=new style 1] (3) at (-1, -1) {};
		\node [style=new style 1] (4) at (1, -1) {};
		\node [style=new style 0] (5) at (-2, -2) {};
		\node [style=new style 0] (6) at (2, -2) {};
	\end{pgfonlayer}
	\begin{pgfonlayer}{edgelayer}
		\draw [style=new edge style 0] (0) to (1);
		\draw [style=new edge style 0] (1) to (2);
		\draw [style=new edge style 0] (2) to (3);
		\draw [style=new edge style 0] (2) to (4);
		\draw [style=new edge style 0] (3) to (5);
		\draw [style=new edge style 0] (4) to (6);
	\end{pgfonlayer}
\end{tikzpicture}
\label{subclaw}
\caption{The subdivided claw}
\end{center}
\end{figure}
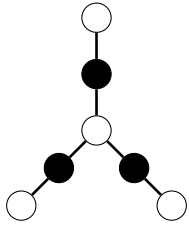

As we are no longer dealing with square patterns, we use bipartite graphs, rather than digraphs to describe the patterns. 
Let $P$ be an $m\times n$ pattern with term-rank $m$. Then 
there exists an $m$-matching, $M$, of $\mathcal{B}(P)$ (that is, $M$ is
a collection of $m$ vertex disjoint edges of $\mathcal{B}(P)$. 
Note that each row vertex is incident to 
(exactly) one edge in $M$, and that there may be more than one $m$-matching of $\mathcal{B}(P)$.
Up to re-labeling of vertices
(or permuting rows and columns of $P$),
we may assume that $P$ is a direct sum of 
matrices whose bigraphs correspond to 
the connected components of $\mathcal{B}(P)$.
Each of these bigraphs has a matching that covers each of its row vertices.

 Throughout filled, respectively unfilled vertices, are row, respectively column vertices. 
The bigraph of the matrix in \Cref{forbiddensun} is
given in 
\Cref{row}. This graph is known as the \textit{subdivided} claw, as it is obtained by subdividing each edge of the claw graph. 
Our first order of business is to describe the bigraph corresponding to a pattern in standard form whose digraph is a weak path.  
A bipartite graph $\mathcal{B}$ is a \textit{Fiedler graph} provided $\mathcal{B}$  consists of 
\begin{itemize}
\item[(a)]  a path, $\gamma$,
$v_1$---$v_2$---$\cdots$--- $v_\ell$
with designated vertices
$1<v_{i_1}<v_{i_2}< \cdots <v_{i_k}< \ell$,
where each of the pairs 
\begin{align*} 
&\mbox{$i_t$ and $i_{t+1}$ ($t=1,\ldots, k-1$)};\\
&\mbox{$1$ and $i_1$}; \mbox{ and }\\
&\mbox{$i_k$  and $\ell$}
\end{align*}
has opposite parity, and  
\item[(b)]
pendant edges $u_i$---$v_i$, where $u_1$, 
\ldots, $u_k$ are the distinct vertices of $G$
not  $\gamma$. 
\item[(c)] $v_1$ and $v_{\ell}$ are pendant vertices of $\mathcal{B}$. 
\end{itemize}
The pendant edges of $\mathcal{B}$ not on $\gamma$ are called \textit{legs}
of $\mathcal{B}$. 
\begin{figure}
\begin{center}

\resizebox{4 in}{!}{
\begin{tikzpicture}
	\begin{pgfonlayer}{nodelayer}
		\node [style=new style 0, fill=black,label=below:$1$] (1) at (1, 0) {};
        \node [style=new style 0, label=below:$1'$] (2) at (2, 0) {};
        \node [style=new style 0,fill=black,label=below:$2$] (3) at (3, 0) {};
        \node [style=new style 0, label=below:$2'$] (4) at (4, 0) {};
        \node [style=new style 0,fill=black,label=below:$3$] (5) at (5, 0) {};
        \node [style=new style 0, label=above:$3'$] (6) at (6, 0) {};
        \node [style=new style 0,fill=black,label=below:$4$] (7) at (7, 0) {};
        \node [style=new style 0,label=below:$4'$] (8) at (8, 0) {};
        \node [style=new style 0,fill=black,label=below:$5$] (9) at (9, 0) {};
        \node [style=new style 0,label=above:$5'$] (10) at (10, 0) {};
         \node [style=new style 0,fill=black, label=below:$6$] (11) at (11, 0) {};
          \node [style=new style 0, label=below:$6'$] (12) at (12, 0) {};
          \node [style=new style 0,label=above:$7'$] (a) at (3, 1) {};
          \node [style=new style 0,fill=black, label=below:$7$] (b) at (6, -1) {};
          \node [style=new style 0, label=above:$8'$] (c) at (7, 1) {};
          \node [style=new style 0,fill=black,label=below:$8$] (d) at (10, -1) {};
	\end{pgfonlayer}
	\begin{pgfonlayer}{edgelayer}
		\draw (1) to (2);
        \draw  (2) to (3);
        \draw  (3) to (4);
        \draw  (4) to (5);
        \draw  (5) to (6);
        \draw  (6) to (7);
        \draw  (7) to (8);
        \draw  (8) to (9);
        \draw  (9) to (10);
        \draw  (10) to (11);
         \draw  (11) to (12);
          \draw (b) to (6);
          \draw (c) to (7);
          \draw (d) to (10);
          \draw (a) to (3);
	\end{pgfonlayer}
    \end{tikzpicture}
}
\\ [6pt]
    \caption{An example of a Fiedler graph}
    \label{fiedgraph}
\end{center}
\end{figure}
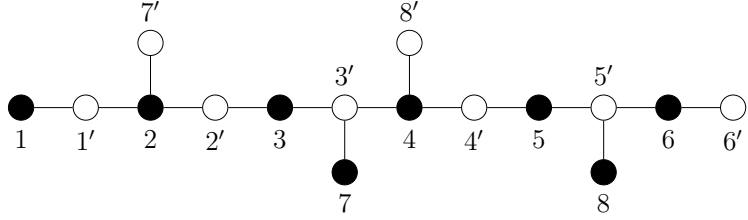
For the example in \Cref{fiedgraph}, 
  $\gamma$ is the path $1$---$1'$---$2$---$2'$---$3$---$3'$---$4$---$4'$---$5$---$5'$---$6$---$6'$;
    the designated vertices are $2$, $3'$, $4$ and $5'$; and 
    the legs are the edges, $2$---$7'$, $3'$---$7$, $4$---$8'$, and $5'$---$8$.
    
Let $P$ be an $m\times m$ pattern of a weak path in standard form. 
Note that if one deletes each row (other than the first or last) 
with exactly one 1 in it, and each column (other than the first or last)  with exactly one 1
in it, one arrives at a matrix whose bigraph 
is a path $\gamma$ on an even number of vertices. Hence the bigraph of $P$ consists of 
$\gamma$ along with pendant edges attached to 
different vertices of $\gamma$. 
As $\mathcal{B}(P)$ has an equal number of row vertices and column vertices, and has a unique perfect matching $M$, every pendant edge is in $M$. Removing the 
edges (and vertices) of these pendant edges results in a disjoint union of subpaths that has a perfect matching.  Hence each of these subpaths of $\gamma$ 
has a perfect matching, and thus has an even number of vertices.  
 In other words, the bigraph of $P$ is a Fiedler graph.  Evidently, the converse also holds.  Thus, the translation of  \Cref{fiedlersquare} into bigraph language is the following. 

\begin{observation}
\label{fiedlergraphdistinct}
Let $P$ by an $m\times m$ matrix of term-rank $m$ whose bigraph is connected.  Then no matrix in $\mathcal{Z}(P)$ has a multiple singular value if and only if $\mathcal{B}(P)$
is a Fiedler graph. 
\end{observation}

We now extend the notion of a Fiedler graph to the setting where there may not be an equal number of rows and columns.
Let $m\leq n$. An $m\times n$ \textit{Fiedler graph} is a graph obtained
from a Fiedler graph $\mathcal{B}$ with $m$ filled vertices and $m$ un-filled vertices
by inserting $n-m$ pendant edges that independently join a new column vertex to a row vertex that is not incident to a leg of $\mathcal{B}$.

\Cref{fiedgraph} is an example of an $8\times 15$ Fiedler graph.
It is obtained from the $8\times 8$ Fiedler graph in \Cref{fiedgraph}
by inserting edges $9'$---$1$, $10'$---$3$, $11'$---$3$, $12'$---$3$, $13'$---$5$,
$14'$---$6$, and $14'$---$6$. 

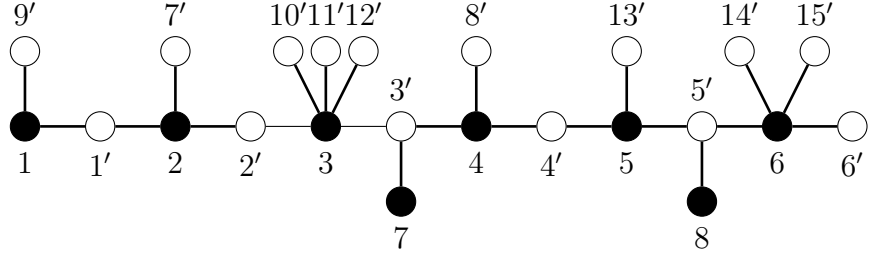
\begin{figure}
\begin{center}
\begin{tikzpicture}[scale=.5]
	\begin{pgfonlayer}{nodelayer}
		\node [style=new style 0, label=below:$2'$] (2) at (-3, 0) {};
		\node [style=new style 1,label=below:$3$] (3) at (-1, 0) {};
		\node [style=new style 0,label=above:$3'$] (4) at (1, 0) {};
		\node [style=new style 1,label=below:$4$] (5) at (3, 0) {};
		\node [style=new style 0, label=below:$4'$] (6) at (5, 0) {};
		\node [style=new style 1, label=below:$1$] (20) at (-9, 0) {};
		\node [style=new style 0, label=below:$1'$] (21) at (-7, 0) {};
		\node [style=new style 1,label=below:$2$] (22) at (-5, 0) {};
		\node [style=new style 1,label=below:$5$] (23) at (7, 0) {};
		\node [style=new style 0,label=above:$5'$] (24) at (9, 0) {};
		\node [style=new style 1, label=below:$6$] (25) at (11, 0) {};
		\node [style=new style 0,label=below:$6'$] (26) at (13, 0) {};
		\node [style=new style 0,label=above:$7'$] (27) at (-5, 2) {};
		\node [style=new style 1, label=below:$7$] (28) at (1, -2) {};
		\node [style=new style 1,label=below:$8$] (29) at (9, -2) {};
		\node [style=new style 0,label=above:$8'$] (30) at (3, 2) {};
		\node [style=new style 0,label=above:$10'$] (31) at (-2, 2) {};
		\node [style=new style 0,label=above:$11'$] (32) at (-1, 2) {};
		\node [style=new style 0,label=above:$12'$] (33) at (0, 2) {};
		\node [style=new style 0,label=above:$15'$] (37) at (12, 2) {};
		\node [style=new style 0,label=above:$14'$] (38) at (10, 2) {};
		\node [style=new style 0,label=above:$13'$] (39) at (7, 2) {};
        \node [style=new style 0, label=above:$9'$] (40) at (-9, 2) {};
	\end{pgfonlayer}
	\begin{pgfonlayer}{edgelayer}
		\draw [style=new edge style 0] (6) to (5);
		\draw [style=new edge style 0] (5) to (4);
		\draw (2) to (3);
		\draw (3) to (4);
		\draw [style=new edge style 0] (20) to (21);
		\draw [style=new edge style 0] (21) to (22);
		\draw [style=new edge style 0] (22) to (2);
		\draw [style=new edge style 0] (6) to (23);
		\draw [style=new edge style 0] (23) to (24);
		\draw [style=new edge style 0] (24) to (25);
		\draw [style=new edge style 0] (25) to (26);
		\draw [style=new edge style 0] (22) to (27);
		\draw [style=new edge style 0] (4) to (28);
		\draw [style=new edge style 0] (5) to (30);
		\draw [style=new edge style 0] (24) to (29);
		\draw [style=new edge style 0] (3) to (31);
		\draw [style=new edge style 0] (3) to (32);
		\draw [style=new edge style 0] (3) to (33);
		\draw [style=new edge style 0] (23) to (39);
		\draw [style=new edge style 0] (25) to (38);
		\draw [style=new edge style 0] (25) to (37);
        \draw [style=new edge style 0] (20) to (40);
	\end{pgfonlayer}
\end{tikzpicture}
\caption{
An $8\times 14$ Fiedler graph.}
\label{exfiedgraph}
\end{center}
\end{figure}
\begin{observation}
\label{fform}
\rm
    If $P$ is a pattern whose bigraph is a $m\times n$ Fiedler graph, then up to permutation of rows and columns we may take $P$ to have the form 
    \[ P = 
    \left[ 
    \begin{array}{c|c}
        F & S  \\
    \end{array}
    \right],
    \]
    where $F$ is the pattern of an $m\times m$ Fiedler graph in standard form,  and $S$ consists of columns of support $1$, corresponding to pendant column vertices. In the \Cref{exfiedgraph}, these pendant vertices correspond to  vertices $9'$, \ldots, $15'$.
\end{observation}
\begin{lemma}
\label{caterpillarsdistinct}
Let $P$ be an $m \times n$ pattern with $m \leq n$ with full term-rank whose bigraph $\mathcal{P}$ is a $m\times n$ Fiedler graph. Let $A \in \mathcal{Z}(P)$. Then the singular values of $A$ are distinct. 
\end{lemma}
\begin{proof}
Without loss of generality, $P$ has the form in \Cref{fform}.
Thus, $A$ has the form  
\[ A = 
    \left[ 
    \begin{array}{c|c}
        B & C  \\
    \end{array}
    \right],
    \]
where $B \in \mathcal{Z}(F)$ and $C\in \mathcal{Z}(S)$.  
Note that the rows of $C$ are mutually orthogonal, and  when 
$i$ corresponds
to a pendant row vertex of $\mathcal{P}$ not on the path $\gamma$ of $\mathcal{B}(F)$, 
then the $i$-th row of $C$ is a row of zeros. Thus $CC\trans$ is a diagonal matrix whose $i$th diagonal 
is $0$ whenever row $i$ is pendant row vertex of $\mathcal{P}$ not on $\gamma$. 

Since $ AA\trans=  
    BB\trans + CC\trans$, the matrices
    $AA\trans$ and $BB\trans$ agree except on possibly $(i,i)$-entries
    where $i$ is not a pendant row vertex.  Now $AA\trans(n,1)$
    is a block lower triangular matrix, each of whose blocks is of the form 
    $((AA^T)[\{j\},\{j+1\}]$  or $((AA^T)[\{j,j+1\}, \{j+1,j+2\}]$ where $j+1$ corresponds to a pendant row vertex.  Hence $BB\trans(n,1)$ 
    is block lower triangular. and its diagonal blocks are identical to those
    of $AA\trans(n,1)$. 
    
    Now let $\sigma$ be a singular value of $A$.  
    Since $F$ has a unique perfect matching, $B$ is invertible. Hence $\sigma>0$. For  $M=AA\trans -\sigma^2 I$,
    $M(n,1)$ is block lower triangular and each diagonal block is a nonzero 
    $1\times 1$ block, or is a block of the form
    \[ 
    \begin{bmatrix}
        a_{i-1,i}a_{i,i}  & a_{i-1,i}a_{i+1,1} \\
        a_{i,i}^2 -\sigma^2 & a_{i,i}a_{i+1,i}
    \end{bmatrix},
    \]
    whose determinant is $\sigma^2>0$.
    Hence $M(n,1)$ is invertible, $M-\sigma^2 I$  has nullity at most $1$, 
    and $\sigma$ is a singular value of $A$ of multiplicity $1$. 
\end{proof}

One of the difficulties with $m\times n$ patterns of term-rank $m$ is that their bigraph
may have more than one $m$-matching.  The observation below describes one way of transitioning from an $m$-matching to 
another $m$-matching. Let $M$ be a matching of a bipartite graph $G$. An \textit{$M$-alternating path}
in $G$ is a path $i_1$---$i_2$--- $\cdots$ ---$i_{2\ell}$ in $G$, 
$i_1$ is not incident to an edge in $M$,  and  $i_{2k}$---$i_{2k+1}$ 
for $k=1,2,\ldots, \ell-1$ are the edges of the path  in $M$.

\begin{observation}
\label{alternating}
Let $G$ be a bipartite graph, $v$ be a vertex of $G$, and  $M$ be a matching of $G\setminus \{v\}$.
Suppose that 
\[ 
\mbox{ $i_1$---$i_2$--- $\cdots$ ---$i_{2\ell}$}
\]
is an $M$-alternating path in $G$ with $v=i_1$.
Then replacing the edges $i_{2k}$---$i_{2k+1}$ 
for $k=1,2,\ldots, \ell-1$ of $M$ by the edges 
$i_{2k-1}$---$i_{2k}$ $k=1,\ldots, \ell-1$
results in a matching of $G\setminus \{i_{2\ell}\}$ with the same number of edges as $M$. 
\end{observation}

\begin{theorem}
\label{bcaterpillars}
Let $P$ be an $m\times n$ pattern with $m\leq n$, term-rank $m$, 
and no columns of all zeros. If $P$ does not allow a matrix with a multiple singular value, then the bigraph of $P$ 
is an $m\times n$ Fiedler graph. 
\end{theorem}
\begin{proof}
   Assume no matrix in $\mathcal{Z}(P)$ 
  has a multiple singular value.  We show that $\mathcal{B}(P)$ is an $m\times n$
  Fiedler graph by establishing a sequence of structural properties of 
  $\mathcal{B}(P)$.  We indicate the end of each proof of a claim by a $\diamond$ .

\bigskip\noindent
{\bf Claim 1.} $\mathcal{B}(P)$ is connected. 
\\[12pt]
    Proof. Suppose not.  Then, after row and column permutations, 
    $P$ is  the direct sum of two nonzero patterns.  There exist matrices with these patterns with the same largest singular values. The resulting direct sum has pattern $P$ and a multiple singular value, contrary to assumption. \hfill $\diamond$ 

\bigskip\noindent
{\bf Claim 2.} Let $\alpha \subseteq \{1,\ldots, n\}$ such that 
$|\alpha|=m$ and the bipartite subgraph $B_{\alpha}$ of $P$ consisting of all row vertices 
and the column vertices in $\alpha$ has an $m$-matching.
Then each connected component of  $B_{\alpha}$
is a Fiedler graph.
\\ [12pt]
Proof. Without loss of generality, we may assume that $\alpha=\{1,\ldots, m\}$.
As $P[\{1,\ldots,m\}]$ has term-rank $m$, each of the connected components of $B_{\alpha}$ 
has a perfect matching.  If one of these components is not a Fiedler graph, then, by   \Cref{fiedlergraphdistinct} and the Direct Sum theorem, there is 
a matrix in $P[\{1,\ldots,m\}]$  with the SSVP and a multiple singular 
value. Hence by the Superpattern theorem, we are led to the contradiction that $P$ allows a matrix with a multiple singular value.
Therefore, each connected component of $B_{\alpha}$ is a Fiedler graph. \hfill $\diamond$ 

\bigskip\noindent
{\bf Claim 3}.
There exists $\alpha \subseteq \{1,\ldots, n\}$ such that 
$|\alpha|=m$ and $B_{\alpha}$ is a connected graph with a perfect matching. 
\\ [12pt]
Proof. Suppose that $B_{\alpha}$ has a perfect matching $M$ and two or more components.
We show how to find an $\alpha'$ such that $B_{\alpha}'$ has a perfect matching 
and fewer connected components. 

Without loss of generality, we may take $\alpha=\{1,\ldots, m\}$.
Claim 1 implies that there exist two connected components $B_1$ and $B_2$ of $B_{\alpha}$
and a column vertex $c$ with $c>m$ such that $c$ has a neighbor $r_1$ in $B_1$
and a neighbor $r_2$ in $B_2$.
Necessarily, there exist vertices $d_1$ and $d_2$ of $B_{\alpha}$
such that {$r_1$---$d_1$} and {$r_2$---$d_2$} are edges in $M$.

If $d_1$ is a pendant vertex of $B_{\alpha}$, then $B_{\alpha \cup \{ c\} \setminus \{ d_1\}}$ has an  $m$-matching
(namely that obtained from $M$ by replacing $r_1$---$d_1$ by $c$---$r_1$), 
and one fewer connected component.
Similarly, we are done if $d_2$ is a pendant vertex in $B_{\alpha}$. 

Now assume that neither $d_1$ nor $d_2$ is a pendant vertex in $B_{\alpha}$. If $r_1$ has degree two or more in $B_\alpha$, 
then $B_{\alpha\cup \{c\}}$ contains a subgraph that is the 
disjoint union of  a subdivided claw with $r_1$ as the center, and a matching 
of size $m-3$. The Direct Sum and Superpattern theorems now 
lead to the contradiction that $P$ allows a matrix with a multiple singular value.  Thus $r_1$, and similarly, $r_2$ are pendant vertices
in $B_{\alpha}$.  As $B_{\alpha \cup \{c\} }$ has an $m$-matching, 
there exists an $M$-alternating path starting at $c$ and containing 
$r_1$. Let $c'$ be the terminal vertex of such a 
path. By \Cref{alternating}, $B_{\{\alpha \cup \{ c\} \setminus{c'} }$
has an $m$-matching and fewer connected components than $B_{\alpha}$.
\hfill $\diamond$

\bigskip 
By Claims 1--3, we can assume that $B_{\{ 1,2,\ldots, m\}}$
is an $m\times m$ Fiedler graph. In particular, it is a tree
with a perfect matching $M$. Let $c>m$ be a column vertex. 

\bigskip\noindent 
{\bf Claim 4. } $c$ is a pendant vertex in $\mathcal{B}(P)$.
\\[12pt]
Proof. Suppose to the contrary that $c$ is not pendant. 
Let $r_1$ and $r_2$ be two of its neighbors.
Let $c_1$ and $c_2$ be the vertices such that $r_1$---$c_1$ and $r_2$---$r_2$' are edges of $M$. 

There exists an $M$-alternating path 
in $B_{\{1,2,\ldots, m,c\}}$ that starts at $c$, contains $r_1$ 
and ends at some pendant column vertex $c'$.  By \Cref{alternating}
$B_{\{1, \ldots, m\}} \cup \{c\} \setminus \{c'\}\}$ contains 
an $M$-matching, but also a cycle. Thus  $B_{\{1, \ldots, m\} \cup \{c\} \setminus \{c'\}\}}$ is not a Fiedler graph, contrary to Claim 1.
Therefore, $c$ is  a pendant vertex in $\mathcal{B}(P)$. \hfill $\diamond$

\bigskip\noindent
Let $\gamma$ be the path of $B_{\{1,\ldots, m\}}$
given in the definition of an ($m\times m$) Fiedler graph.

\bigskip\noindent
{\bf Claim 5.}  The neighbor $r$  of $c$ is not a pendant vertex
of $B_{\{1,2,\ldots, m\}}$, unless $r$ is the 
first or last vertex of the path $\gamma$. 
\\ [12pt]
Proof. 
Suppose not. Let $c'$ be the neighbor of $r$ in $B_{\{1,\ldots, n\}}$.  As $\mathcal{B}$ has a perfect matching $M$,
and $c'$ is not pendant,  $r$--$c'$ is in $M$, $c'$ has two additional neighbors $r_1$ and $r_2$ on $\gamma$  Let $r_1$---$c_1$ and $r_2$---$c_2$ be the edges of $M$ that contain $r_1$ and $r_2$.
Note that  the induced subgraph $H$ with vertices $c,r,c',r_1,r_2,c_1,c_2$
is a subdivided claw. Hence $B_{\{1,\ldots, m,c\}}$
has a subgraph that is the disjoint union of a subdivided claw 
and a $(m-3)$-matching.  As before, this contradicts that $P$ does not
allow a matrix with a multiple singular value. \hfill $\diamond$

\bigskip\noindent
Claims 1-5 now imply that $\mathcal{B}(P)$ is an $m\times n$ Fiedler graph.
\end{proof}

\Cref{caterpillarsdistinct} and \Cref{bcaterpillars} now give a complete characterization of $m\times n$ patterns of term-rank $n$ that do not 
allow a matrix with a multiple singular value. 
\begin{corollary}
  Let $P$ be an $m\times n$ pattern with term-rank $m$ and no column of zeros. 
  Then no matrix with pattern $P$ has a multiple singular value if and only if the bipartite graph of $P$ is an $m\times n$ Fiedler graph.  
\end{corollary}
\section{Square Patterns with less than Full Term-Rank}

Let $S$ be an $m \times n$ pattern. If the term-rank of $S$ is $m-2$ or less, then every matrix $A\in \Z(S)$ has $0$ as a singular value with multiplicity 2. Hence, in the context of requiring distinct singular values, we can focus our attention on $m\times n$ patterns with $m\leq n$ and  term-rank $m-1$ or $m$.  The case of term-rank $m$ is resolved in the previous sections. 

In this section, we initiate the study of the case where 
$S$ has term-rank $m-1$.  
Such $S$ do not have the SSVP, which precludes the use of the Superpattern theorem and the Direct Sum theorem.  Consequently, one cannot use the SSVP to find a list of forbidden subgraphs to classify the patterns without full term-rank that require distinct singular values.  Despite the limitations with handling such patterns, the term-rank constraint does impose some exploitable structure.

We prove two results that may be of further use in characterizing such patterns that require all simple singular values.  The first, \Cref{n-2form}, 
characterizes the patterns that do not allow 
0 as a multiple singular value. The second, \Cref{orthogborder} establishes a way to construct families of such patterns $S$ that allow a multiple nonzero singular values. 

We begin with some needed definitions.
A \textit{line} of a matrix $A$ is a  column or a row of $A$. A \textit{line cover} of $A$ is a set of lines of $A$ that contain all nonzero entries of $A$. A consequence of the K\"onig-Egerv\'ary theorem is that an $m\times n$ matrix $A$ has term-rank $t$ if and only if 
the minimum number of lines in a cover of 
$A$ is $t$. We will use $m_A(\sigma)$ to denote the \textit{multiplicity} of $\sigma$ as a singular value of $A$. The following characterizes the  $m \times n$ patterns with $m\leq n$ and term-rank do not allow a matrix with $0$ as a singular value. This is essentially an immediate consequence of Theorem 3.4 of \cite{HershkowitzSchneider}.  We include a proof here for the reader's convenience. 
\begin{lemma}
\label{requireli}
Let $P$ be an $m\times n$ pattern. Then 
every matrix with pattern $P$ has rank  $m$ if and only 
if the rows and columns of $P$ can be permuted to obtain 
a matrix of the form 
\[ 
\begin{bmatrix}
\begin{array}{c|c}
S & T
\end{array}
\end{bmatrix},
\]
where $S$ is an $m\times m$ upper triangular matrix 
with all diagonal entries equal to $1$.
\end{lemma}

\begin{proof}
We prove the forward direction by induction on $m$.
A $1 \times n$ pattern $P$ requires rank 1 if and only if it is nonzero. This is equivalent to $P$ having the desired form. 
Assume $m>1$ and proceed by induction. Let $P$ be an $m\times n$
pattern that requires rank $m$.

If each column of $P$ has either zero or at least two nonzeros, 
then there is a matrix with pattern $P$ whose column sums are 0, 
contradicting the rank $m$ requirement. Thus, some column of $P$ has exactly one nonzero entry.
Without  loss of generality, we may assume that 
\[
P= \begin{bmatrix}
\begin{array}{c|c}
    1 & \mathbf u\trans \\ \hline 
    \mathbf 0 & Q
    \end{array} 
\end{bmatrix}
\]
for some $(m-1)\times (n-1)$ matrix $Q$. 
Since $P$ requires rank $m$, $Q$ requires rank $m-1$.
By the inductive assumption, the rows and columns of $Q$ 
can be permuted so that the leading $(m-1)\times (m-1)$ 
submatrix of $Q$ is upper triangular with all diagonal entries $1$.
It follows that $P$'s rows and columns can be permuted so that 
its leading $m\times m$ submatrix is upper triangular with all diagonal entries equal to $1$.

The converse direction follows by noting 
every matrix with pattern $S$ has linearly independent rows. 
\end{proof}

The next result characterizes $m\times n $
zero-nonzero patterns with $m\leq n$ that require 
rank $m-1$.  This is essentially an immediate consequence of Theorem 3.9 of \cite{HershkowitzSchneider}.

\begin{lemma}
\label{n-2form}
    Let $P$ be an $m \times n$ pattern with $m \leq n$ of that requires rank $m-1$. Then 
    there exist positive integers $r$ and $s$ with $r+s=n+1$ such that up to row and column permutations  $P$ has the form 
    \[ 
    \widehat{P}=
    \begin{bmatrix}
    \begin{array}{c|c}
        P_{11} & O \\ \hline 
        P_{12} & P_{22} \end{array}
    \end{bmatrix},\]
    where 
    \begin{itemize}
        \item[\rm (a)]
    $P_{11}$ is $r \times (r-1)$ with term-rank $r-1$
    and its leading $(r-1)\times (r-1)$ submatrix
    is upper triangular with all diagonal entries equal to $1$,  and
    \item [\rm (b)]
    $P_{22}$ is $(m-r) \times s$ with term-rank $m-r$ and its leading $(m-r)\times (m-r)$ submatrix is 
    upper triangular with all diagonal entries equal to $1$.
    \end{itemize}
    Furthermore, if 
    \[
    M = \left\lbrack 
    \begin{array}{c|c}
         A & O  \\
         \hline
         B & C 
    \end{array}
    \right\rbrack \in \mathcal{Z}(\widehat{P}),
    \]
then  
    the columns of $A$ are linearly independent, and 
    the rows of $C$ are linearly independent.
\end{lemma}
\begin{proof}
Since $P$ has term-rank $m-1$, there is a line cover of $P$ with $m-1$ lines. Hence, up to row and column permutations, 
 $P$ has the form 
 \[
    P = \left\lbrack 
    \begin{array}{c|c}
         P_{11} & O  \\
         \hline
         P_{21} & P_{22} 
    \end{array}
    \right\rbrack,
    \]
    where $P_{11}$ is $r \times (r-1)$. Hence,  $P_{21}$ is $(m-r) \times (r-1)$, and $P_{22}$ is  $(m-r) \times (n-m+1)$. Since $P$ requires rank $m-1$, $P_{11}$ requires rank $r-1$.
    Hence by \Cref{requireli}, statement (a) holds. 

    Let  $C$ be a matrix with pattern $P_{22}$ and  let $\mathbf y\trans$ be a vector 
    with $\mathbf y{\trans} C=\mathbf{0} \trans$. Let $B$ be any  matrix with pattern $P_{21}$ and $A$ 
    any matrix 
    with pattern $P_{11}$. Since (a) holds, $A$ has rank $r-1$. Thus there exists a vector $\mathbf z \trans$ such that $\mathbf z \trans A=-\mathbf y\trans B$.  Therefore, 
    \[ 
    [ \mathbf{z}\trans  \, \,  \mathbf{y}\trans] 
    \begin{bmatrix}
        A & O \\
        B & C 
    \end{bmatrix}
    = \mathbf 0\trans.
    \]
    Since $A$ has more rows than columns, there is a nonzero vector 
    $\mathbf{w}\trans$ such that $\mathbf {w} \trans A=\mathbf 0$.
    Hence 
    \[
     [ \mathbf{w}\trans  \, \,  \mathbf{0}\trans] 
    \begin{bmatrix}
        A & O \\
        B & C 
    \end{bmatrix}
    = \mathbf 0\trans.
    \]
    As 
    \[ \begin{bmatrix} 
        A & O \\ 
        B & C 
        \end{bmatrix}
    \]
    has pattern $P$, its nullity is at most one. This implies that 
    $\mathbf y\trans= \mathbf 0$. Hence every matrix with pattern $P_{22}$ has linearly independent rows.   By \Cref{requireli},
    (b) holds. 

    The furthermore statement readily follows from (a) and (b). 
\end{proof}

\begin{Example}\rm
    Let $A$ be a matrix of the form
    \[ A = 
    \left\lbrack
    \begin{array}{c|cc}
         a & 0 & 0  \\
         b & 0 & 0\\
         \hline
         c & d & e
    \end{array}
    \right\rbrack,
    \]
    where each of $a$, $b$, $c$ and $d$ is nonzero.
 Set
    \[ Q_1 =  
    \frac{1}{\sqrt{a^2+b^2}}\left\lbrack
    \begin{array}{rr}
         b& -a\\
         a & b
    \end{array} 
    \right\rbrack\oplus [1],
    \]
    and
    \[ Q_2 =  
    [1] \oplus \frac{1}{\sqrt{d^2+e^2}}\left\lbrack
    \begin{array}{rr}
       d & e\\
       e & -d
    \end{array}
    \right\rbrack.
    \] 
Then 
    \[
    Q_1AQ_2 = \left\lbrack
\begin{array}{ccc}
     0 & 0 & 0   \\
     \frac{a^2}{\sqrt{a^2+b^2}} & 0 & 0   \\
     c & \frac{d^2}{\sqrt{d^2+e^2}} & 0   \\
\end{array}
    \right\rbrack.
    \]
Since both $Q_1$ and $Q_2$ are orthogonal, $Q_1AQ_2$ has the same singular values as $A$. The matrix $Q_1AQ_2$  has a multiple singular value if and only if
\[
(Q_1AQ_2)(1,3)= \left\lbrack
\begin{array}{rr}
  \frac{a^2}{\sqrt{a^2+b^2}} & 0    \\
     c & \frac{d^2}{\sqrt{d^2+e^2}} 
\end{array}
\right\rbrack,
\]
has a multiple singular value. The latter occurs  if and only if  the rows of $(Q_1AQ_2)(1,3)$ are orthogonal 
and have the same length.

Consequently, the pattern 
\[ 
\begin{bmatrix}
    1 & 0 & 0 \\
    1 & 0 & 0 \\
    1 & 1 & 1 
\end{bmatrix}
\]
has term-rank $2$ and requires distinct singular values. 
Also, the pattern 
\[ 
\begin{bmatrix}
    1 & 0 & 0 \\
    1 & 0 & 0 \\
    0 & 1 & 1 
\end{bmatrix}
\]
has term-rank $2$ and allows a matrix with a multiple singular value. $\diamond$
\end{Example}

More generally, we have the following result. 

\begin{lemma}

\label{borderow}
    Let $A$ be an $m \times n $  matrix with $m \leq n$  and a singular value $\sigma$ such that $m_A(\sigma) = k$. Let $\ell < k$ and let $\widehat{A}$ be the $(m+\ell)\times n$ matrix formed by bordering $A$ below by an $\ell \times n$ matrix $U$. Then $m_{\widehat{A}}(\sigma) \geq k-\ell$. 
\end{lemma}
\begin{proof}
This follows from the Cauchy interlacing inequalities and the observation that 
$AA\trans$ is a principal submatrix of 
$\widehat{A} \widehat{A}\trans$.
\end{proof}

\begin{theorem}
\label{orthogborder}
    Let $Q$ be an $n\times n$ matrix with the SSVP whose rows are mutually orthogonal and of equal Euclidean length $\ell$.
    Let $P$ be the zero-nonzero pattern of $Q$. 
    Let $k< n$, $p \geq k$, and $S$ be a superpattern of the pattern
    \[
    \left[
    \begin{array}{r|l}
         P & O  \\
         \hline
        O & O_{k,p}
    \end{array}
    \right].
    \]
    Then there exists a matrix $A$ with zero-nonzero pattern $S$ 
    with  $m_{A}(\ell) \geq n-k$.
\end{theorem}

\begin{proof}
    Since the rows  of  $Q$ are nonzero and mutually orthogonal, its rows are linearly independent.
    Thus, by Theorem 3.6 of \cite{cheung2025strongsingularvalueproperty},  the matrix $\left[ \begin{array}{c|c}Q  &  O_{n,p} \end{array} \right] $ has the SSVP, and $\Sigma(Q) = \Sigma(\left[ \begin{array}{c|c}Q  &  O_{n,p} \end{array} \right])$. Note 
    that $S[\{ 1, \ldots, n\}, \{1, \ldots, n+p\}]$  is a superpattern of $\left[ \begin{array}{c|c}P  &  O_{n,p} \end{array} \right] $. 
    By the Superpattern theorem, there exists a $B$ with zero-nonzero pattern $S[\{ 1, \ldots, n\}, \{1, \ldots, n+p\}]$  
   with  $m_B(\ell) = n$. Now, by \cref{borderow},bordering $B$ below by any $k\times (n + p)$ matrix results in a matrix $A$ with pattern $S$ and $m_A(\ell) \geq n - k$. 
\end{proof}

Using \cref{orthogborder} we can generate families  of $n \times n$ patterns with term-rank $n-1$ that have 1 as a singular value of multiplicity at least 2. 
A \textit{full lower Hessenberg pattern},
$H_n$ is the $n\times n$ pattern whose entries are 0 above the first superdiagonal and whose other entries are all 1.  By Corollary 2.3 of (\cite{Cheon1999}), $H_n$ is the pattern 
of an orthogonal matrix.

\begin{lemma}
\label{hessvp}
    Let $n \geq 2$. Then every matrix whose pattern is $H_n$ has the SSVP.
\end{lemma}
\begin{proof}
    Let $X$ be an $n \times n$ matrix such that $A\trans X$ and $X A\trans$ are symmetric and $A \circ X = O$. The third condition implies that $X$ is strictly upper triangular. Hence $A\trans X$ is strictly upper triangular. We show that the $j$-th superdiagonal, that is, the entries of $A\trans X$ in positions $(i,i+j)$) of $X$ are zero by induction on $j$. Since $A \circ X = O$, the first superdiagonal of $X$ is all 0's by definition. Now, suppose $j > 1$ and the  $k$th superdiagonal of $X$ is zero for all $k < j$. Then the entries of the $j$-th superdiagonal of $A\trans X$ are of the form $a_{i,i+j}x_{i,i+j}$. Since $A\trans X$ is symmetric, and $A\trans X$ is strictly upper triangular, $a_{i,i+1}x_{i,i+j} = 0$. However, every $a_{i,i+1} \neq 0$ for $j \geq 1$, thus $x_{i,i+j} = 0$. Hence, $X =O$ and $A$ has the SSVP. 
    \end{proof}

\begin{Example}
\rm 
    Consider the $8\times 8$ pattern
    \[ P = 
    \left[ 
    \begin{array}{cccccccc}
        1 & 0 & 0 & 0 & 0 & 0 & 0 & 0  \\
        1 & 1 & 0 & 0 & 0 & 0 & 0 & 0  \\
        1 & 1 & 1 & 0 & 0 & 0 & 0 & 0 \\
        1 & 1 & 1 & 1 & 0 & 0 & 0 & 0 \\
        1 & 1 & 1 & 1 & 1 & 0 & 0 & 0 \\
        1 & 1 & 1 & 1 & 1 & 0 & 0 & 0 \\
        1 & 1 & 1 & 1 & 1 & 1 & 0 & 1\\
        1 & 1 & 1 & 1 & 1 & 1 & 1 & 1\\
    \end{array}
    \right].
    \]
  Note that $P\lbrack \{2,3,4,5,6\},\{1,2,3,4,5\}\rbrack = H_5$. Since $H_5$ allows orthogonal rows, there is a matrix $A$ whose pattern is $H_5$ with $m_A(1) = 5$.  
  By \cref{hessvp} $A$ has the SSVP, and by \cref{orthogborder} there is a matrix $M$ whose pattern is $P$ with $m_M(1)  \geq 2$. Note that $M$ has term-rank  $7$. Since this pattern is not full term-rank, it does not have the SSVP. So such a construction cannot be be obtained using either the Direct Sum theorem or the Matrix Liberation theorem.
\hfill $\diamond$  \end{Example}

More generally, we can find entire families of examples of $n \times n$ square matrices with $1$ as singular value of multiplicity at least 2 and term-rank $n-1$ using the following construction. Let $H_n^+$ denote the $(n+1) \times n$ pattern formed by stacking the standard basis vector $\mathbf{e}\trans_1$ on top of  $H_n$.

\begin{Example}
\rm
Let $n \geq 4$, $S$ be a superpattern of $O$, and let $R$ be an $(n-3) \times (n-2)$ pattern whose last column is arbitrary, and whose leading $(n-3) \times (n-3)$ subpattern has $1$'s on its diagonal, $0$'s above its first diagonal, and all other entries arbitrary. Then, the pattern,
    \[ P = 
\left[ 
\begin{array}{c|c}
    H^+_{n} & O\\
    \hline
    S& R
\end{array}
\right]
    \]
    has a realization $M$ which has $1$ as a singular value of multiplicity at least $2$, and term-rank $2n-2$. 

   As mentioned previously, there is a matrix $A$ with pattern $H_n$ whose rows are orthogonal, so $m_A(1) = n$. By \cref{hessvp}, $A$ also has the SSVP, so the desired result follows from \cref{orthogborder}
\hfill $\diamond$  \end{Example}

\bibliographystyle{acm} 
\bibliography{Refs.bib}
\end{document}